\documentclass[11pt,êeqno]{article}
\usepackage{amsfonts}
\usepackage[tbtags]{amsmath}
\topskip=\baselineskip%
\newcommand{\nfty}{N \to \infty}
\newcommand{\lr}{\left(}
\newcommand{\lb}{\left[}
\newcommand{\rb}{\right]}
\newcommand{\rr}{\right)}
\newcommand{\lf}{\left\{}
\newcommand{\rf}{\right\}}
\newcommand{\lv}{\left|}
\newcommand{\rv}{\right|}
\newcommand{\lp}{\left.}
\newcommand{\rp}{\right.}
\newcommand{\bel}{\begin{lemma}}
\newcommand{\beql}[1]{\begin{equation}\label{#1}}
\newcommand{\eel}{\end{lemma}}
\newcommand{\xkn}{X_{k:N}}
\newcommand{\xmn}{X_{m:N}}
\newcommand{\ukn}{U_{k:N}}
\newcommand{\umn}{U_{m:N}}
\newcommand{\xin}{X_{i:N}}
\newcommand{\uin}{U_{i:N}}
\newcommand{\xknd}{X^2_{k:N}}
\newcommand{\xmnd}{X^2_{m:N}}

\newcommand{\xind}{X^2_{i:N}}

\newcommand{\xia}{\xi_{\alpha}}
\newcommand{\xib}{\xi_{\beta}}
\newcommand{\lon}{\log N/N}

\newcommand{\na}{N_{\alpha}}
\newcommand{\nau}{N_{\alpha, u}}
\newcommand{\nb}{N_{\beta}}

\newcommand{\win}{W_{i:N}}
\newcommand{\wi}{W_i}
\newcommand{\ha}{\frac 1{f(\xi_{\alpha})}}
\newcommand{\hb}{\frac 1{f(\xi_{\beta})}}
\newcommand{\alp}{\alpha}
\newcommand{\be}{\beta}
\newcommand{\si}{\sigma}
\newcommand{\Ha}{\xi_{\alpha}}
\newcommand{\Hb}{\xi_{\beta}}
\newcommand{\Har}{\xi_{\alpha}^r}
\newcommand{\Hbr}{\xi_{\beta}^r}
\newcommand{\Harl}{\xi_{\alpha}^{r-1}}
\newcommand{\Hbrl}{\xi_{\beta}^{r-1}}
\newcommand{\Had}{\xi_{\alpha}^2}
\newcommand{\Hbd}{\xi_{\beta}^2}
\newcommand{\dd}{\partial}
\newtheorem{lemma}{Lemma}

\title{\huge \sf The Empirical Edgeworth Expansion \\for a Studentized Trimmed Mean
}
\author{{\sf Nadezhda Gribkova$^{1}$} \ and  \ \ {\sf Roelof Helmers$^{2}$}}
\date{ 
{\small \it
\centerline{$^{1}$~St.Petersburg~State~University,~Mathematics~and~Mechanics~Faculty,}
\centerline{198504,~St.Petersburg,~Stary~Peterhof,~Universitetsky~pr.~28,~Russia}
\centerline{E-mail:~nv.gribkova@gmail.com}
\centerline{$^{2}$ Center~for~Mathematics~and~Computer~Science}
\centerline{P.O.Box~94079,~1090~GB~Amsterdam,~The~Netherlands}
\centerline{E-mail:~helmers@cwi.nl}
}} 
\begin{document}
\maketitle
\begin{abstract}
{\small \it We establish the validity of the empirical Edgeworth expansion
(EE) for a Studentized trimmed mean, under the sole condition that
the underlying distribution function of the observations satisfies
a local smoothness condition near the two quantiles where the
trimming occurs. A simple explicit formula for the $N^{-1/2}$ term
(correcting for skewness and bias; $N$ being the sample size) of
the EE will be given.
%In particular our result supplements
%previous work by Hall and Padmanabhan \cite{halp} and Putter
%and van Zwet [12].
The proof is based on a U-statistic type approximation and also
uses a version of Bahadur's \cite{bah} representation for sample
quantiles. }
\end{abstract}

\noindent 1991 Mathematics Subject Classification: Primary 62E20, 62G30;
Secondary 60F05.\\[3mm]
{\em Keywords and phrases}: Empirical Edgeworth expansions, Studentized trimmed means.\\

\subsection*{1. Introduction}

\ \ \ \ \ \ The trimmed mean is a well known estimator of a
location parameter. Its asymptotic properties were studied by many
authors (see \cite{bick}, \cite{bile}, \cite{bj}, \cite{halp},
\cite{hjqz}, \cite{stig}, \cite{tukey}, and the references therein).
The main reason for applying the trimmed mean is
robustness (see \cite{hamp}, \cite{huber}). The limit distribution
of the trimmed mean for an arbitrary population distribution was
found by Stigler \cite{stig}. Specifically he has shown that in
order for the trimmed mean to be asymptotically normal, it is
necessary and sufficient that the sample is trimmed at sample
quantiles for which the corresponding population quantiles are
uniquely defined \cite{stig}.

In this paper we study the second-order asymptotic properties of
the distribution of the trimmed mean, as well as of the Studentized trimmed
mean in view of its practical relevance (construction of confidence intervals, hypothesis tests
etc.).

We establish the validity of the empirical Edgeworth expansion
(EEE) for a Studentized trimmed mean, under the sole condition
that underlying distribution function ($df$) of the observations
satisfies a local smoothness condition near the two quantiles
where the trimming occurs. In particular our result supplements
previous work by Hall and Padmanabhan \cite{halp} and Putter and
van Zwet \cite{putz}. The existence of an Edgeworth expansion (EE)
for a Studentized trimmed mean was also obtained by Hall and
Padmanabhan \cite{halp}, but these authors wrote that the "first
term in an Edgeworth expansion is very complex and so it will not
be written down explicitly". They suggested to replace analytical
difficulties by bootstrap simulation. In contrast, in the present
paper we show that our method of proof gives a simple explicit
formula for the $N^{-1/2}$- term (correcting for skewness and
bias; $N$ being the sample size) of the Edgeworth expansion.

The proof of our result is based on a $U-$statistic type
approximation (cf. also Bickel et al \cite{bigz}, Helmers
 \cite{helm1}-\cite{helm2}, Putter and van Zwet \cite{putz})
 and also uses a version of Bahadur's
\cite{bah} representation for sample quantiles. Our $U-$statistic
type approximation is slightly different from the one given by the
first two terms of the Hoeffding decomposition and approximates
the trimmed mean with a remainder of the classical Bahadur's order
$N^{-3/4}\log N^{5/4}$ (cf. (4.4)-(4.5), Sect.4). The first order
linear term of our $U-$statistic approximation is a sum of
independent identically distributed (i.i.d.) Winsorized random
variables. The structure of the quadratic term of the second order
is connected with a Bahadur type property of the order statistics
close to the sample quantile (cf. lemma 3.2, Sect.3). We will also
show (cf. Lemma A.2, Appendix) that our result cannot be obtained
as a consequence of a general result on Edgeworth expansions for
Studentized symmetric statistics (Theorem 1.2, \cite{putz}) of
Putter and van Zwet.

In Section 2, we formulate and discuss our main results on EE and
EEE. In Section 3, we state and prove Bahadur's type lemmas. Next,
in Section 4, we construct $U-$statistic type approximation for
the trimmed mean and prove the result on EE for the normalized
trimmed mean. In Section 5, the corresponding stochastic
approximation for a plug-in estimator, which is used to construct
a Studentized trimmed mean is established, and the result on EE
for a Studentized trimmed mean is proved. Finally, in Section 6,
we prove some lemmas on the consistency of our estimators of the
unknown parameters appearing in the formula of one-term EE and
establish a rate of convergence. In the Appendix, we establish an
asymptotic approximation for the bias of trimmed mean in
estimating of the corresponding location parameter, and prove that
our results on EE and EEE for a Studentized trimmed mean can not
be inferred from results of Putter and van Zwet \cite{putz} for
Studentized symmetric statistics.

\subsection*{2. The main results}
\ \ \ \  Let $X_1,\dots ,X_N$ be i.i.d. real-valued random
variables (r.v.) with common $df$ $F$, and let \
$X_{1:N}\le\cdots\le X_{N:N}$ denote the corresponding order
statistics. Consider the trimmed mean given by
$$
T_N=\frac 1{([\be N]-[\alp N])}\sum_{i=[\alp N]+1}^{[\be N]}\xin\
, \leqno (2.1)
$$
where $0<\alp<\be<1$ are any fixed numbers and [$\cdot$]
represents the greatest integer function. Let $F^{-1}(u)= \inf \{
x: F(x) \ge u \}$, $0<u< 1$, denote the left-continuous inverse
function of $df$ $F$ and put $\frac d{du}F^{-1}(u)=1/f(F^{-1}(u))$
to be its derivative, when the density $f=F'$ exists and
$f(F^{-1}(u))>0$. Let
$$
\xi_{\nu}=F^{-1}(\nu),
$$
$0<\nu<1$, be the $\nu$-th quantile of $F$. Define a function
$$
Q(u)=\left\{
\begin{array}{cl}
\Ha \ ,& u\le \alp ,\\
F^{-1}(u)\ ,& \alp < u\le \be ,\\
\Hb \ ,& \be < u .
\end{array}
\right.
%\leqno (1.3)
$$
Let $\wi$, $i=1,\dots , N$, denote $X_i$ Winsorized outside
of $(\Ha,\Hb]$, that is
$$
W_i=\left\{
\begin{array}{cl}
\Ha \ ,& X_i\le \Ha, \\
X_i \ ,& \Ha < X_i \le \Hb ,\\
\Hb \ ,& \Hb < X_i .
\end{array}
\right. \leqno (2.2)
$$
Then $\wi \stackrel{\rm d}{=} Q(U_i)$, $i=1,\dots , N$, where
$U_i$ are independent r.v.'s with uniform $(0,1)$ distribution.
Define
$$
\mu_W=\int_0^1Q(u)\, du, \quad  \si_W^2= \int_0^1(Q(u)-\mu_W)^2 \,
du, \quad \gamma_{3,W}=\int_0^1(Q(u)-\mu_W)^3 \, du . \leqno (2.3)
$$
Put
$$
\delta_{2,W}=-\alp^2\ha[\mu_W-\Ha]^2+(1-\be)^2 \hb [\mu_W-\Hb]^2 .
\leqno (2.4)
$$
Suppose that $\xia\neq\xib$ (that is $\xia$ is not an atom with
mass at least $(\be - \alp)$ for the distribution $F$), then the
$W_i$'s are not degenerate. Define  real numbers $\lambda_1$ and
$\lambda_2$ by
$$
\lambda_1= {\gamma_{3,W}}/{\si_W^3}, \quad \quad
\lambda_2={\delta_{2,W}}/{\si_W^3}. \leqno (2.5)
$$
We need no moment assumptions about the distribution $F$
and to normalize $T_N$ we use
$$
\mu (\alp,\be)=\frac 1{\be -\alp} \int_{\alp}^{\be}F^{-1}(u)\, du
\leqno (2.6)
$$
as a location parameter and $(\be -\alp)^{-1}\si_W$ (the root of
the asymptotic variance, cf.(4.8)) as a scale parameter. Note that
$T_N$ often serves as a statistical estimator for the parameter
$\mu (\alp,\be)$, the population trimmed mean.

Now we show why moments are not needed. Take some fixed $\Delta
>0$ and define auxiliary i.i.d. Winsorized r.v.'s $X'_i=\max
(\Ha-\Delta,\min(X_i,\Hb+\Delta))$. Let $\xin'$, $i=1,\dots,N$,
denote the corresponding order statistics. Introduce an auxiliary
trimmed mean $T'_N=\frac 1{([\be N]-[\alp N])}\sum_{i=[\alp
N]+1}^{[\be N]}\xin'$, and note that
$$
T_N =T'_N \quad \mbox{if} \quad
\{X_{[\alp N]+1:N} \ge \Ha -\Delta\} \cap \{X_{[\be N]:N} \le \Hb +\Delta\} .
$$
If $F$ has a positive and continuous density in neighborhoods of
$\Ha$ and $\Hb$, then, by Bernstein's inequality $P(X_{[\alp
N]+1:N}<\Ha -\Delta) + P(X_{[\be N]:N}>\Hb +\Delta)=O(exp(-cN))$,
as $\nfty$, where $c>0$ is constant independent of $N$ . Therefore
$$
\sup_{x\in R} |P(T_N\le x) - P(T'_N\le x)|= O(e^{-cN}) \leqno
(2.7)
$$
and when proving our results we can replace with impunity $T_N$ by $T'_N$,
which has finite moments of the arbitrary order.

In absence of any moment assumptions, our formulas for the
$N^{-1/2}$ term of the Edgeworth expansions contains a bias term.
Define the quantity
$$
\be_N=\frac 1N \lf- (\alp N-[\alp N])\lr \frac {}{}\mu (\alp ,
\be) -\Ha\rr- \frac 12 \alp(1-\alp)\ha \rp \leqno (2.8)
$$
$$
\qquad +\lp  (\be N-[\be N])\lr \frac {}{}\mu (\alp , \be) -\Hb
\rr+ \frac 12 \be (1-\be)\hb\rf .
$$
Note that when both $\alp N$ and $\be N$ are integer valued, the
bias term has a very simple form: $\be_N =\frac 1{2N} \lf -\frac
{\alp(1-\alp)}{f(\xia)} + \frac {\be(1-\be)}{f(\xib)}\rf$.
Moreover, in case $\alp = 1-\be$ and $f(\xia)=f(\xib)$ (when the
distribution $F$ is symmetric, for example), the bias term
vanishes.

We show (cf. Lemma A.1, Appendix) that if the conditions of our
Theorem 2.1 are satisfied, then for an arbitrary $\Delta >0$
$$
b_N= (\be -\alp ) (ET'_N - \mu (\alp ,\be)) = \be_N + O(N^{-3/2})
\leqno (2.9)
$$
as $\nfty$.(cf. (2.7)) Note also that the bias term (2.8) does not
depend on the auxiliary quantity $\Delta$.

Define
$$
F_{T_N}(x) = P\lr \frac{N^{1/2}(T_N - \mu (\alp , \be))}{(\be -
\alp)^{-1}\si_W} \le x \rr \leqno (2.10)
$$
to be the distribution function of the normalized trimmed mean.
Using the notation of Putter and van Zwet \cite{putz}, we shall
show that the Edgeworth expansion for $F_{T_N}(x)$ is given by
$$
G_N(x)= \Phi (x) - \frac {\phi (x)}{6\sqrt{N}}\lr (\lambda_1+3
\lambda_2) (x^2-1) +6 N \frac {\be_N}{\si_W}\rr, \leqno (2.11)
$$
where $\Phi$ is the standard normal distribution function, $\phi
=\Phi'$. The quantity $(\lambda_1+3 \lambda_2)N^{-1/2}$ serves as
an approximation to the third cumulant of $\frac{N^{1/2}(T'_N -
\mu (\alp , \be))} {(\be - \alp)^{-1}\si_W}$, moreover $\lambda_1
N^{-1/2}$ is the approximation to the third cumulant of the
$L_2$-projection of the normalized trimmed mean, which close to
$N^{-1/2}\si_W^{-1} \sum_1^N \wi$ - a sum of N i.i.d. Winsorized
r.v.'s (cf.~Sect.4, below), and $3 \lambda_2 N^{-1/2}$ is due to
the $U-$statistic type approximation to $T_N$.

Here is our first result: an Edgeworth expansion for a normalized trimmed mean.
%\begin{theorem}

\bigskip
{\sc Theorem} 2.1.\ \ {\it Suppose that $f=F'$ exists in
neighborhoods of the points $\xi_{\alp}$ and $\xi_{\be}$ and
satisfies a Lipschitz condition. In addition we assume that
$f(\xi_{\nu})>0$, $\nu=\alp,\be$. Then
$$
\sup_{x\in R}|F_{T_N}(x) - G_N(x)|= o(N^{-1/2}), \leqno (2.12)
$$
as $\nfty$.}
%\end{theorem}

\bigskip
Theorem 2.1 can be viewed as a version of the Edgeworth expansion
for the trimmed mean obtained by Bjerve \cite{bj} in his
unpublished Berkeley Ph.D. thesis (cf. also Helmers \cite{helm1}).
Our method of proof is completely different from Bjerve's, as he
used a conditioning argument to reduce a trimmed mean to a sum of
i.i.d. r.v.'s, conditionally given the values of $X_{[\alp
N]+1:N}$ and $X_{[\be N]:N}$, while in contrast we essentially
show that $T_N$ can be approximated by a $U-$statistic $U_N$; the
remainder $T_N-U_N$ can be shown to be of negligible order for our
purposes by an application of a version of Bahadur \cite{bah}
representation for sample quantiles.
%It should be noted that the another way to
%prove (1.12) is to check the condition (1.15) of Putter and van
%Zwet [12] and to apply their Thm.1.1. We can show that their
%condition(1.15) is satisfied for a normalized trimmed mean,
%provided we strengthen the local smoothness condition used in our
%Thm 1.1 slightly.

Next we state our result on the validity of one-term Edgeworth
expansion for the Studentized trimmed mean. Define plug in estimators for $\mu_W$
and $\si^2_W$ by
$$
\hat \mu_W =\frac kN \xkn +\frac 1N \sum_{i=k+1}^{m-1}\xin + \frac
{N-m+1}{N}\xmn , \leqno (2.13)
$$
and
$$
S_N^2=\lr \frac kN \xknd +\frac 1N \sum_{i=k+1}^{m-1}\xind + \frac
{N-m+1}N \xmnd \rr - \hat \mu^2_W  \leqno (2.14)
$$
with $k=[\alp N]+1$ and $m=[\be N]$. Let
$$
F_{N,S}(x) =P\lr \frac{N^{1/2}(T_N - \mu (\alp , \be))}{(\be -
\alp)^{-1}S_N} \le x \rr \leqno (2.15)
$$
denote the $df$ of a Studentized trimmed mean. Define
$$
H_N(x)= \Phi (x) + \frac {\phi (x)}{6\sqrt{N}}\lr \frac {}{}
(2x^2+1)\lambda_1+3(x^2+1)\lambda_2 -6 N \frac {\be_N}{\si_W}\rr .
\leqno (2.16)
$$
Our main result is:
%\begin{theorem}

\bigskip
{\sc Theorem} 2.2. \ \ {\it Suppose that the conditions of Theorem
2.1 are satisfied. Then
$$
\sup_{x\in R}|F_{N,S}(x)-H_N(x)|= o(N^{-1/2}), \leqno (2.17)
$$
as $\nfty$ .
}
%\end{theorem}
%\begin{remark}

\bigskip

As already indicated in our introduction the existence of an
Edgeworth expansion for $F_{N,S}$ was proved by Hall and
Padmanabhan \cite{halp}. In (2.16) and (2.17) we give the precise
and simple explicit form of the Edgeworth expansion for $F_{N,S}$.
In fact formally the form of our $H_N$ (cf.(2.16)) coincides with
the one given on p.1545 of Putter and van Zwet \cite{putz}.
However, our Theorem 2.2 can not be inferred from the result of
Putter and van Zwet \cite{putz}: the second condition in
assumption (1.18) of Putter and van Zwet [16, p.1542], is not
satisfied for our $T_N$, that is, for a Studentized trimmed mean
(cf.~Lemma A.2, Appendix). Our conjecture is that also the  first
condition in their assumption (1.18) is not satisfied , but this
seems rather difficult to check for a Studentized trimmed mean.

\bigskip
{\sc Remark} 2.1. \ \ It is clear from the proofs of Theorems 2.1
and 2.2 that the order of the remainder term which we really
obtain in relations (2.12) and (2.17) is $O((\log
N)^{5/4}/N^{3/4})$, as $\nfty$.
%\end{remark}

\bigskip
To obtain empirical Edgeworth expansions (cf. Helmers
\cite{helm2}, Putter and van Zwet \cite{putz}) we replace
$\lambda_1$, $\lambda_2$, $\be_N$ and $\si_W$ in (2.11) and (2.16)
by statistical estimates. The estimation of $\lambda_1$ is
straightforward. Let us define
$$
\hat \lambda_1 = S_N^{-3} \hat \gamma_{3,W}  \qquad \qquad \qquad
\qquad \qquad \qquad \qquad
$$
$$
=S_N^{-3}\lr \frac kN (\xkn -\hat \mu_W)^3 +\frac 1N
\sum_{i=k+1}^{m-1} (\xin -\hat \mu_W)^3+ \frac {N-m+1}N (\xmn
-\hat \mu_W)^3 \rr
$$ ($\hat \mu_W$ and $S_N$ were defined in
(2.13) and (2.14)) to be an estimate for $\lambda_1$. As to
$\lambda_2$ and $\be_N$, we first have to estimate the values of
density $f(\xia)$ and $f(\xib)$. We shall use kernel estimators
with a simple step-like kernel. Put $g(x)=I_{\{|x|\le 1/2\}}$.
Take the width of kernel $\delta = N^{-1/4}$ and put
$g_{\delta}(x)=\frac 1{\delta} g\lr \frac x{\delta}\rr = \frac
1{\delta}I_{\{|x|\le \delta /2\}}$, where $\int_{-\infty}^{\infty}
g_{\delta}(x) \, dx = 1$. Then our estimates for values of density
at the quantiles where trimming occurs will be the following:
$$
\hat f (\xi_{\nu})=\frac 1N \sum_{i=1}^N g_{\delta}(X_i -X_{r:N})
= N^{-3/4} \sum_{i=1}^N I_{\{ 2N^{1/4}|X_i -X_{r:N} |\le 1 \}},
\leqno (2.18)
$$
where $\nu=\alp$ and $r=k$ or $\nu=\be$ and $r=m$ respectively.
Our estimates of $f(\xia)$ and $f(\xib)$ are rather simple ones
and sufficient for our purposes (cf. also Reiss [17, p.262]). One
easily obtains the following estimates for $\lambda_2$ and
$\be_N$:
$$
\hat \lambda_2=S_N^{-3}\lf -\alp^2(\hat f(\xia))^{-1} [\hat \mu_W -\xkn]^2 +
(1-\be)^2 (\hat f(\xib))^{-1} [\hat \mu_W - \xmn ]^2 \rf ,
$$
$$
\hat \be_N = \frac 1{N}\lf -(\alp N - [\alp N])\lr \frac {}{} T_N
- \xkn \rr - \frac 12 \alp (1-\alp) (\hat f(\xia))^{-1}  \rp
$$
$$
\qquad + \lp  (\be N - [\be N]) \lr \frac {}{} T_N - \xmn \rr
+\frac 12 \be (1-\be) (\hat f(\xia))^{-1} \rf .
$$
When the conditions of Theorem 2.1 are satisfied, the estimates
$\hat \lambda_1$, $\hat \lambda_2$ and $\hat \be_N$ are consistent
estimators of the corresponding quantities $\lambda_1$,
$\lambda_2$ and $\be_N$ (cf.~Sect.6). Replacing the latter
quantities by these estimates in formulas (2.11) and (2.16), we
obtain the empirical Edgeworth expansions:
$$
\hat G_N(x)= \Phi (x) - \frac {\phi (x)}{6\sqrt{N}}\lr (\hat \lambda_1+
3 \hat \lambda_2 )
(x^2-1) + 6 N \frac {\hat \be_N}{S_N} \rr,
$$
$$
\hat H_N(x)= \Phi (x) + \frac {\phi (x)}{6\sqrt{N}}\lr \frac {}{}
(2x^2+1)\hat \lambda_1+3(x^2+1)\hat \lambda_2 -
6 N\frac {\hat \be_N}{S_N}\rr .
$$

Our result, establishing the validity of the empirical Edgeworth
expansions, is given by the following assertion.
%\begin{theorem}

\bigskip
{\sc Theorem} 2.3. \ \ {\it Suppose that the conditions of Theorem
2.1 hold. Then
$$
\sup_{x\in R}|F_{T_N} (x)-\hat G_N(x)| = o_p \lr \frac 1{\sqrt
{N}}\rr, \leqno (2.19)
$$
$$
\sup_{x\in R}|F_{N,S}(x) -\hat H_N(x)| = o_p \lr \frac 1{\sqrt
{N}}\rr. \leqno (2.20)
$$
as $\nfty$.
}
%\end{theorem}

\bigskip

{\sc Remark} 2.2. \ \ It is clear from Remark 2.1 and the Lemma's
6.1 and 6.2 that we can strengthen (2.19) and (2.20) to
$\sup_{x\in R}|F_{T_N} (x)-\hat G_N(x)|= O\lr (\log
N)^{5/4}N^{-3/4}\rr$ with probability $1- O\lr N^{-c}\rr$, for
every $c>0$, as $\nfty$, and similarly, $\sup_{x\in
R}|F_{N,S}(x)-\hat{H}_N(x)|= O\lr (\log N)^{5/4}N^{-3/4}\rr$,
except on a set with probability $O(N^{-c})$, for every
$c>0$.\\[3mm]
 To conclude this section we remark that an alternative
way of approximating $F_{T_N}$ or $F_{N,S}$ accurately is to use
saddlepoint methods. In Helmers et al \cite{hjqz} saddlepoint
approximations were established rigorously for the trimmed mean
and the Studentized trimmed mean. Compared with the Edgeworth
expansions derived in the present paper, the saddlepoint
approximations will typically behave better in the far tail of the
distribution. An advantage of empirical Edgeworth expansions is
that they are much easier to compute.

%\end{remark}
\subsection*{3. Auxiliary results}
%{\bf 2. Auxiliary results.}
\ \ \ \ \  Define the binomial r.v. $ \na= \sharp \{i : X_i \le
\xia \}$ , where \  $0<\alp <1. $

The following lemma is a version of Bahadur's \cite{bah}
representation (cf. also Theorem 6.3.1, Reiss \cite{reis}) for the
sample quantile. In this section $k$ denotes an integer satisfying
$k=\alp N +O(1)$, $\nfty $.

\bigskip
{\sc Lemma} 3.1. \ \ {\it Suppose that $f=F'$ exists and is
positive and Lipschitz in neighborhood of $\xia$. Let $G$ be a
function defined in a neighborhood of $\xia$ and $g=G'$ exists and
satisfies a Lipschitz condition. Then
$$
G(\xkn) = G(\xia) - \frac {\na - \alpha N}{N} g(\xia)/f(\xia ) +
R_N, \leqno (3.1)
$$
where
$$
P(|R_N|>A(\lon )^{3/4}) =O(N^{-c}),  \leqno (3.2)
$$
as $\nfty$, for every $c>0$ and some $A>0$, not depending on $N$.}
%\eel

\bigskip

We omit the proof because the lemma is essentially known and its
proof requires similar arguments,which will also be used in the
proof of Lemma 3.2. Our proof of the next lemma will use the
following fact: conditional on $\na$ the order statistics
$X_{1:N},\dots ,X_{\na:N}$ are distributed as  $\na$ i.i.d. r.v.'s
with distribution function $F(x)/\alp$, \ $x\le \xia$. Though this
fact is more or less  known, we add
a brief explanation of it.
Let $U_1,\dots,U_N$ are independent
r.v.'s uniformly distributed on $(0,1)$ and $U_{1,N},\dots
,U_{N,N}$ denote the corresponding order statistics. Put $\nau=
\sharp \{i : U_i \le \alp \}$. Since $\xin \stackrel{\rm d}{=}
F^{-1}(\uin)$ and $\na \stackrel{\rm d}{=} \nau\,$, it is enough
to prove the assertion for the uniform distribution.
\ \ First consider the case $\nau= N$. Take
arbitrary $0<u_1\le \cdots \le u_N < \alp$ and write
\begin{equation*}
\label{ap_1}
\begin{split}
P(U_{1:N} \le u_1,\dots,U_{\nau:N}\le u_N \mid\nau =N)
=\frac{P(U_{1:N} \le u_1,\dots,U_{n:N}\le u_N)}{\alp^N}\\
=\frac{N!}{\alp^N}\int_0^{u_1}\int_{u_1}^{u_2}\dots
\int_{u_{N-1}}^{u_N} d\,x_1 d\,x_2\dots d\,x_N,
\end{split}
\end{equation*}
and the latter is $d.f.$ of the order statistics corresponding to
the sample of \ $N$ \ independent $(0,\alp)$-uniform distributed
r.v.'s. \ \  Now let  $\nau=k<N$ and  $F_{i,N}(u)=P(U_{i:N}\le u)$ be a~$df$ of $i$-th order statistic, put $P_N(k)=P(\nau=k)={N\choose k} \alp^k(1-\alp)^{N-k}$. Then we can write
\begin{equation*}
\label{ap_2}
%\begin{split}
P(U_{1:N} \le u_1,\dots,U_{\nau :N}\le u_k \mid \nau =k)
= \frac {P(U_{1:N} \le u_1,\dots,U_{k :N}\le u_k, U_{k+1:N}>\alp)}{P_N(k)}.
\end{equation*}
The probability in the nominator on the r.h.s. of  the latter formula is equal to
$
\int_{\alp}^1 \!\!\! P\big(U_{1:N} \le
u_1,\dots,U_{k:N} \le u_k \mid U_{k+1:N}=v\big)\,d F_{k+1,N}(v),
$
and by the Markov property of order statistics the latter quantity equals
\begin{equation*}
\label{ap_3}
\begin{split}
&\phantom{ =}\int_{\alp}^1 \lr \frac{k!}{v^k} \int_0^{u_1}\int_{u_1}^{u_2}\dots \int_{u_{k-1}}^{u_k} d\,x_1
d\,x_2\dots d\,x_k\rr \,d F_{k+1,N}(v)\\
&= \frac{k!}{\alp^k}\lr \int_0^{u_1}\int_{u_1}^{u_2}\dots \int_{u_{k-1}}^{u_k} d\,x_1
d\,x_2\dots d\,x_k \rr\times \alp^k\int_{\alp}^1 \frac 1{v^k} \,d F_{k+1,N}(v),
\end{split}
\end{equation*}
and since $\alp^k\int_{\alp}^1 \frac 1{v^k} \,d F_{k+1,N}(v)=\alp^k\int_{\alp}^1 \frac{(1-v)^{N-k-1}}{B(k+1,N-k)} \,d v
= {N\choose k}\alp^k (1-\alp)^{N-k}=P_N(k)$, where $B(k+1,N-k)=k!(N-k-1)!/N!$, we obtain that conditional probability we consider is equal to
\begin{equation*}
\frac{k!}{\alp^k}\int_0^{u_1}\int_{u_1}^{u_2}\dots
\int_{u_{k-1}}^{u_k} d\,x_1 d\,x_2\dots d\,x_k,
\end{equation*}
which  corresponds to the $(0,\alp)$-uniform distribution.

To state next lemma we shall adopt the following notation. Let
$\sum_{i=k}^m (.)_i = sign[m-k]
\sum_{i= k\wedge m}^{k\vee m} (.)_i $
for all integer $k$ and $m$.
%\bel

\bigskip
{\sc Lemma} 3.2. \ \ {\it Suppose that the conditions of lemma 3.1
are satisfied. Then
$$
\frac 1N \sum_{i=k}^{\na } \lr G(\xin )-G(\xia ) \rr = - \frac
{(\na - \alp N)^2}{2N^2} g(\xia )/f(\xia ) + R_N, \leqno (3.3)
$$
where
$$
P(|R_N|>A (\lon )^{5/4})= O(N^{-c}),\leqno (3.4)
$$
as $\nfty$ for every $c>0$ with some $A>0$, not depending on $N$.}
%\eel

\bigskip
This lemma extends and sharpens the relations (3.2) and (3.3)
given (for the case $G(x)=x$) in Hall and Padmanabhan \cite{halp}.
Note also that the factor $(1-\alp)^{-1}$ in formula (3.2) and
$(1-\be)^{-1}$ in formula (3.3) (see Hall and Padmanabhan
\cite{halp}) should be omitted. We  apply this lemma several
times: to approximate the trimmed mean (cf. lemma 4.1), its
asymptotic variance (cf. lemma 5.1) and its asymptotic third
moment (cf. Thm.2.3 and lemma 6.2).

{\sc Corollary} 3.1. {\it Suppose that $f=F'$ exists and is
positive and Lipschitz in a neighborhood of $\xia$. Then
$$
\frac 1N \sum_{i=k}^{\na } (\xin - \xia ) = - \frac {(\na - \alp
N)^2}{2N^2}\ha + R_{N,1},
$$
$$
\frac 1N \sum_{i=k}^{\na } \lr X_{i,N}^2 - \xi_{\alp }^2 \rr = -
\frac {(\na - \alp N)^2}{N^2} \Ha \ha + R_{N,2},
$$
where $R_{N,i}$, $i=1,2$, satisfy} (3.4).

{\sc Proof}. We begin by writing (cf.(3.3))
$$
R_N= \frac 1N \sum_{i=k}^{\na} \lr G(\xin ) - G(\xia )\rr +\frac
{(\na - \alp N)^2}{2N^2} g(\xia )/f(\xia ).\leqno (3.5)
$$
Now we will check that $R_N$ satisfies (3.4). Let, as before,
$U_1,\dots,U_N$ denote independent r.v.'s uniformly distributed on
$(0,1)$ and let $U_{1,N},\dots ,U_{N,N}$ denote the corresponding
order statistics. Since the joint distribution of $\xin$,
$(i=k,\dots,\na)$ and $\na$ coincides with the joint distribution
of $F^{-1}(\uin)$ $(i=k,\dots,\nau)$ and $\nau$, where $\nau =
\sharp \{i : U_i \le \alp \}$, it of course suffices to verify
that
$$
\frac 1N \sum_{i=k}^{\nau} [G(F^{-1}(\uin )) - G(F^{-1}(\alp ))] +
\frac {(\nau - \alp N)^2}{2N^2} g(\xia )/f(\xia ) \leqno (3.6)
$$
satisfies (3.4). By our smoothness condition the first term of
(3.6) equals
$$
\frac 1N g(\xia)/f(\xia ) \sum_{i=k}^{\nau } (\uin - \alp) +
R_{N,3}, \leqno (3.7)
$$
where
$$
|R_{N,3}|\le \frac CN \sum_{i=k \wedge \nau}^{k \vee \nau}(\uin
-\alp)^2 \le \frac {C|k-\nau|}N \lb (\ukn -\alp )^2 \vee
(U_{\nau,N}-\alp)^2\rb \leqno (3.8)
$$
with $C$ is equal to the Lipschitz constant of function
$g(F^{-1}(u))/f(F^{-1}(u))$ (we neglect here the event that $\ukn$
does not belong to the neighborhood of $\alp$ where smoothness
conditions hold, as this probability is of the order
$O(exp(-cN))$, as $N\to\infty$ for some $c>0$, cf. the
introduction). Let us fix an arbitrary $c>0$ and note that
$$
P\lr (\alp -U_{\nau,N})^2>A_1\lon \rr \le P\lr U_{\nau +1,N}-U_{\nau,N}
> (A_1\lon)^{1/2} \rr
$$
$$
= P\lr U_{1:N}> (A_1\lon)^{1/2} \rr = O(N^{-c}). \leqno (3.9)
$$
Here and elsewhere $A_j$ denote the positive constants which do not depend
on $N$. Besides, by Bernstein's inequality
$$
P(|\nau -k|>(A_2 N\log N)^{1/2})=O(N^{-c}), \leqno (3.10)
$$
with $A_2=2c\alp (1-\alp)$, and by lemma 3.1.1, Reiss \cite{reis}
$$
P((\ukn -\alp )^2 > A_3 \lon )= O(N^{-c}),
$$
as $\nfty $. Therefore (3.8) implies that
$$
P(|R_{N,3}|> A_4 (\lon )^{3/2})=O(N^{-c}) \leqno (3.11)
$$
with $A_4=CA_2\max(A_1,A_3)$. Next we consider the dominant term
on the r.h.s. of (3.7). By (3.10) we can bound our quantities on
the event $E=\lf \omega :|\nau - k|<(A_2N\log N)^{1/2} \rf$. Fix
$N$ and $\nau $ for which the event $E$ holds true. Without loss
of generality let $k\le \nau$. Note that conditional on $\nau $
the order statistic $\uin $, $k\le i \le \nau$, is distributed as
$i$-th order statistic of the sample of size $\nau $ from the
uniform on $(0,\alp )$ distribution and $E(\uin | \nau ) = \frac
{\alp i}{\nau +1}$, for $i=k,\dots , \nau$. Write
$$
\frac 1N \sum_{i=k}^{\nau } (\uin -\alp )= \frac 1N \lb
\sum_{i=k}^{\nau } (\uin -\frac {\alp i}{\nau +1})
+\sum_{i=k}^{\nau } (\frac {\alp i}{\nau +1} - \alp) \rb \leqno
(3.12)
$$
$$
= \frac 1N \sum_{i=k}^{\nau } (\uin -\frac {\alp i}{\nau +1}) -
\frac {\alp (\nau -\alp N)^2}{2N\nau } + O((\log
N)^{1/2}N^{-3/2}).
$$
For the second term on the r.h.s. of (3.12) we have
$$
- \frac {\alp (\nau -\alp N)^2}{2N\nau } = - \frac {(\nau -\alp
N)^2}{2N} \frac {\alp }{\alp N+(\nau -\alp N)} =- \frac {(\nau
-\alp N)^2}{2N^2} +R_{N,4}, \leqno (3.13)
$$
where in view of (3.10)
$$
P(|R_{N,4}|>A_5 (\lon )^{3/2})=O(N^{-c}) \leqno (3.14)
$$
as $\nfty$ with $A_5=A_2^3$. For the first term on the r.h.s. of
(3.12) we can write
$$
\frac 1N \lv \sum_{i=k}^{\nau } \lr \uin -\frac {\alp i}{\nau
+1}\rr \rv \le \frac {\nau -k+1}N \max_{k\le i\le \nau} \lv \uin -
\frac {\alp i}{\nau +1}\rv . \leqno (3.15)
$$
Note that we suppose that the event $E$ holds true and (without
loss of generality) that $k\le \nau $ (otherwise a similar
argument with respect to $(\alp,1)$ instead $(0,\alp)$) will do).
Fix an arbitrary $c_1>c+1/2$ and note that conditional on $\nau $
the variance of the order statistic $\uin $, $k\le i \le \nau$, is
equal to $\frac {\alp^2 i (\nau -i+1)}{(\nau +1)^2(\nau +2)} =
O\lr (\log N)^{1/2} N^{-3/2} \rr$. By lemma 3.1.1, Reiss
\cite{reis}, we obtain that uniformly for $k\le i \le \nau $
$$
P\lr \lv \uin - \frac {\alp i}{\nau +1}\rv > A_6 (\lon )^{3/4}|
\nau \rr= O(N^{-c_1}), \leqno (3.16)
$$
as $\nfty$.Relations (3.15) and (3.16) together imply
$$
P\lr \frac 1N \lv \sum_{i=k}^{\nau } (\uin - \frac {\alp i}{\nau +1}) \rv
>(A_2)^{1/2}A_6\lp (\lon )^{5/4} \rv \nau \rr \le \leqno (3.17)
$$
$$
(A_2N\log N)^{1/2}O(N^{-c_1})=O(N^{-c}),
$$
as $\nfty $. Now (3.3) and (3.4) follows from (3.5)--(3.7),
(3.11)--(3.14) and (3.17). The lemma is proved. $\Box$

\subsection*{4. Proof of Theorem 2.1}
%{\bf 3. Proof of Theorem 1.1.}
\ \ \ \ \ \ To begin with let us note that we can replace $T_N$
(cf. (2.1)) by
$$
N^{-1/2}\sum_{i=k}^m \xin,   \leqno (4.1)
$$
where $k= [\alp N]+1$, $m=[\be N]$, $0<\alp <\be < 1$. Note that
though $T_N$ in (4.1) is of different order than in (2.1), this
will affect only the bias term (see Lemma A.1, Appendix), and we
shall take that into account whenever needed. Define
$I_{\nu}(X_i)=I_{\{X_i\le \xi_{\nu}\}}$, where
$\xi_{\nu}=F^{-1}(\nu)$, \ $0<\nu<1$, and $I_A$ is the indicator
of event $A$. Then for the Winsorized r.v. $W_i$ (cf. (2.2)) we
can write
$$
W_i=X_iI_{\be}(X_i)(1-I_{\alp}(X_i))+\Ha I_{\alp}(X_i) +\Hb
(1-I_{\be}(X_i)). \leqno (4.2)
$$
Recall that $\mu_W$, $\si^2_W$, $\gamma_{3,W}$ denote first three
cumulants of r.v. $W_1$ (cf.(2.3)). Define  a $U-$statistic of
degree 2 by
$$
L_N + U_N = \sum_{i=1}^N L_{N,i} + \sum_{1 \le i}\sum_{< j \le N}
U_{N,(i,j)}, \leqno (4.3)
$$
where
$$
L_{N,i} = \frac 1{\sqrt{N}} (\wi - \mu_W) \qquad \qquad \qquad
\qquad \qquad \leqno (4.4)
$$
$$
=\frac 1{\sqrt{N}}\lb \frac {}{} X_iI_{\be}(X_i)(1-I_{\alp}(X_i))+
\Ha I_{\alp}(X_i) +\Hb (1-I_{\be}(X_i)) - \mu_W \rb ,
$$
$$
U_{N,(i,j)} = \frac 1 {N\sqrt{N}} \lb \frac {}{} - \ha (I_{\alp
}(X_i)-\alp )(I_{\alp }(X_j)-\alp )  \rp \leqno (4.5)
$$
$$
\qquad + \lp \frac {}{} \hb (I_{\be }(X_i)-\be )(I_{\be }(X_j)-\be
) \rb .
$$
Note that
$$
EL_{N,i}=0 \leqno (4.6)
$$
for all $i=1,\dots ,N$ and
$$
E U_{N,(i,j)} =0, \qquad E(L_{N,i}U_{N,(i,j)})=0 \leqno (4.7)
$$
for all $i,j=1,\dots ,N$ $(i\ne j)$. Using (4.4)--(4.7), we easily
check that
$$
\si^2_{L_N+U_N} = E(L_N+U_N)^2= E(L_N^2)+O(N^{-1})=\si^2_W +
O(N^{-1}), \leqno (4.8)
$$
and also that
$$
E(L_N+U_N)^3= E(L_N^3) + 3 E(L_N^2 U_N) + O(N^{-3/2})   \qquad
\qquad \ \ \leqno (4.9)
$$
$$
= \frac 1{\sqrt{N}} \gamma_{3,W}+ 3\frac 1{\sqrt{N}} \lf -\ha \lb
\frac {}{} E((W_1-\mu_W) (I_{\alp}(X_1)-\alp))\rb^2 \rp \quad
\qquad
$$
$$
+ \lp \hb \lb \frac {}{} E((W_1-\mu_W) (I_{\be}(X_1)-\be))\rb^2
\rf + O(N^{-3/2})  \quad \qquad \qquad \qquad
$$
$$
\quad =\frac 1{\sqrt{N}} \gamma_{3,W} +3\frac 1{\sqrt{N}} \lb -\ha
\alp^2 \frac {}{} [\Ha -\mu_W]^2 + \hb (1-\be )^2 [\Hb -\mu_W]^2
\rb
$$
$$
+ O(N^{-3/2}).
$$
Relations (4.8) and (4.9) imply that
$$
E\lr \frac {L_N+U_N}{\si_{(L_N+U_N})} \rr ^3 = \frac {\lambda_1 +
3\lambda_2}{\sqrt{N}} + O(N^{-3/2}), \leqno (4.10)
$$
with  $\lambda_1$ and $\lambda_2$ as in (2.5).

The next lemma ensures that the approximation of $T_N$ by a
$U-$statistic of the form (4.3) has a remainder of classical
Bahadur's order of magnitude $N^{-3/4}(\log N)^{5/4}$.

\bigskip
{\sc Lemma} 4.1. \ \ {\it Suppose that the conditions of Theorem
2.1 hold. Then
$$
P\lr |T_N-ET'_N-(L_N+U_N)| > A(\log N)^{5/4}N^{-3/4}\rr =
O(N^{-c}) \leqno (4.11)
$$
as $\nfty$, for every $c>0$ with some $A>0$ independent on $N$.
}

\bigskip
{\sc Proof of Lemma 4.1}. Let $\win $, $i=1,\dots ,N$, denote the
order statistics, corresponding to $W_1,\dots ,W_N$. Put $N_{\nu
}= \sharp \{X_i : X_i \le \xi_{\nu } \}, \ 0<\nu<1$. Then
$$
\win =\left\{
\begin{array}{cl}
\Ha \ ,& i\le \na, \\
\xin \ ,& \na < i \le \nb ,\\
\Hb \ ,& i> \nb .
\end{array}
\right.
$$
Now note that
$$
T_N - \frac 1{\sqrt{N}} \sum_{i=1}^N W_i = \frac 1{\sqrt{N}}\lr
\sum_{i=k}^m \xin - \na \Ha - \sum_{i=\na +1}^{\nb} \xin -
(N-\nb)\Hb \rr
$$
$$
=\frac 1{\sqrt{N}}\lf sign[\na-(k-1)]\sum_{i=k\wedge (\na
+1)}^{\na \vee (k-1)} \xin - sign(\nb - m) \sum_{i=(m\wedge \nb)
+1}^{m\vee \nb} \xin \rp
$$
$$
-\lp \frac {}{} \na \Ha - (N-\nb)\Hb \rf = \frac 1{\sqrt{N}}\lf
sign[\na-(k-1)]\sum_{i=k\wedge (\na +1)}^{\na \vee (k-1)} (\xin
-\Ha)\rp \ \
$$
$$
-\lp  sign(\nb - m) \sum_{i=(m\wedge \nb)+1}^{m\vee \nb} (\xin
-\Hb) -(k-1) \Ha - (N-m)\Hb \rf  \qquad \qquad \ \
$$
$$
=-\frac{(\na-\alp N)^2}{2N\sqrt{N}}\ha + \frac{(\nb-\be
N)^2}{2N\sqrt{N}}\hb - \frac {k-1}{\sqrt{N}}\Ha - \frac
{N-m}{\sqrt{N}}\Hb +R_N, \quad \
$$
where by Lemma 3.2
$$
P\lr |R_N| > A(\log N)^{5/4}N^{-3/4}\rr = O(N^{-c}) \leqno (4.12)
$$
as $\nfty$, for every $c>0$ with some $A>0$ independent of $N$.
Define
$$
Q_N =-\frac{(\na-\alp N)^2}{2N\sqrt{N}}\ha + \frac{(\nb-\be
N)^2}{2N\sqrt{N}}\hb \qquad \qquad \qquad \qquad \quad
$$
$$
\qquad \qquad =\frac 1{2N\sqrt{N}}\lf -\lb \sum_{i=1}^N
(I_{\alp}(X_i)-\alp)\rb^2\ha + \lb \sum_{i=1}^N
(I_{\be}(X_i)-\be)\rb^2\hb \rf .
$$
It is clear that $Q_N$ is a symmetric polynomial of degree two
with
$$
E(Q_N)= \frac 1{2\sqrt{N}}\lf \frac {}{} -\alp (1-\alp)\ha +
\be (1-\be)\hb\rf .
$$
Note that
$$
E \frac 1{\sqrt{N}} \sum_{i=1}^N W_i = \sqrt{N} \mu_W =
\sqrt{N}\lr \frac {}{} (\be-\alp)\mu (\alp,\be)+\alp \Ha +(1-\be)
\Hb \rr . \leqno (4.13)
$$
Next we can write
$$
T_N=L_N+Q_N - EQ_N  +\sqrt{N}\lr \frac {}{} (\be-\alp)\mu
(\alp,\be) +\alp \Ha +(1-\be) \Hb \rr \qquad \quad \leqno (4.14)
$$
$$
\qquad \  - \frac {k-1}{\sqrt {N}}\Ha - \frac {N-m}{\sqrt {N}}\Hb
+ \frac 1{2\sqrt{N}}\lf \frac {}{} -\alp (1-\alp)\ha + \be
(1-\be)\hb\rf +R_N
$$
$$
= L_N+Q_N - EQ_N+\sqrt{N}(\be-\alp)\mu (\alp,\be) +\frac
1{\sqrt{N}} \lf \frac {}{} -(k-1-\alp N)\Ha  \rp  \qquad
$$
$$
- \lp  \frac 12 \ha \alp (1-\alp) +(m-\be N)\Hb + \frac 12 \hb \be
(1-\be) \rf +R_N. \qquad \qquad
$$
Let us compare the expression within curly brackets on the r.h.s.
of (4.14) with the formula for $B_2$ (the bias term for $T'_N=
N^{-1}\sum_{i=k}^m \xin'$) (cf.(A.3), Appendix). As a result we
obtain the following formula:
$$
T_N - \sqrt{N}(\be-\alp)\mu (\alp,\be) = L_N +Q_N-EQ_N
+\sqrt{N}B_2 +R_N, \leqno (4.15)
$$
with $R_N$ as in (4.14) plus $O(N^{-1})$ (cf. (A.3), Appendix).
Note that $R_N$ satisfies (4.12), and as $T'_N$ is normalized by
$N^{-1/2}$ in this lemma, we have $ET'_N- \sqrt{N}(\be-\alp)\mu
(\alp,\be)=\sqrt{N}B_2$ (cf. lemma A.1, Appendix). So, relation
(4.15) implies
$$
T_N - ET'_N = L_N+Q_N-EQ_N  +R_N. \leqno (4.16)
$$
For the quantity $Q_N-EQ_N$ we can write
$$
Q_N - EQ_N = U_N + \frac {\bar r_N}{2\sqrt{N}}, \leqno (4.17)
$$
where $U_N$ as in (4.3) and
$$
\bar r_N = \frac 1N \sum_{i=1}^N \lf \frac {}{} -\ha [(I_{\alp
}(X_i)-\alp )^2-\alp (1-\alp )] \rp
$$
$$
+ \lp \frac {}{}  \hb [(I_{\be }(X_i)-\be )^2-\be (1-\be )] \rf.
$$
Note that $\bar r_N$ is the average of $N$ i.i.d. bounded and
centered ($E \bar r_N=0$) r.v.'s, and by Hoeffding's inequality
\cite{hoef} we find that
$$
P\lr |\bar r_N | >A (\lon )^{1/2} \rr = O(N^{-c}) \leqno (4.18)
$$
for every $c>0$ and some $A>0$, not depending on $N$. Therefore
$\frac 12 \bar r_N /\sqrt{N}$ on the r.h.s. of (4.17) is
negligible for our purposes. Relations (4.12) and (4.16)--(4.18)
together imply (4.11). The lemma is proved.  $\Box$

\bigskip
{\sc Remark} 4.1. \ \ The first linear term of our U-statistic
approximation to $T_N$ is a sum of i.i.d. Winsorized r.v.'s $W_i$.
A simple argument involving formula (2.10) for the
$L_2$-projection (i.e. ,the first term of the Hoeffding
decomposition) given in [16, p.1548], tells us that our leading
term is slightly different from the one given by the  Hoeffding
decomposition. The same fact holds true for the second quadratic
term in our U -statistic approximation to the trimmed mean.

\bigskip
{\sc Proof of Theorem} 2.1. Using the Lemma 4.1 and Lemma A.1
(cf.Appendix), for the $df$ of $T_N$ (cf. (2.1)) defined by (2.10)
we can write
$$
F_{T_N}(x) = P\lf \frac {N^{1/2} (\be - \alp)(T_N-ET'_N)}{\si_W }
\le x - \frac {N^{1/2}\be_N +O(N^{-1})}{\si_W}\rf  \qquad \leqno
(4.19)
$$
$$
\qquad \qquad \ \qquad \quad = P\lf \frac {L_N+U_N}{\si_W} \le
\frac {[\be N]- [\alp N]}{(\be - \alp)N} \lr x -\frac
{N^{1/2}\be_N}{\si_W} + O(N^{-1}) \rr - \frac {R_N}{\si_W}\rf
$$
$$
\qquad \qquad \ \quad \quad = P\lf \frac {L_N+U_N}{\si_W} \le  x
(1+O(N^{-1})) -\frac {N^{1/2}\be_N}{\si_W}-\frac {R_N}{\si_W}
+O(N^{-1}) \rf ,
$$
where $L_N+U_N$ (cf.(4.3)) is $U$-statistic of degree two with the
canonical functions
$$
g_N(x) = E(L_N+U_N |X_1=x) \qquad \qquad \qquad \qquad \qquad
\qquad \qquad \qquad \qquad \qquad
$$
$$
= \frac 1{\sqrt{N}}[xI_{\be}(x)(1-I_{\alp}(x))+ \Ha I_{\alp}(x) +
\Hb (1-I_{\be}(x)) - \mu_W], \qquad
$$
$$
\psi_N(x,y) = E(L_N+U_N |X_1=x, X_2=y) - g_N(x) - g_N(y)  \qquad
\qquad \qquad \qquad \qquad \ \
$$
$$
\qquad \qquad \quad   = \frac 1{N\sqrt{N}} \lb \frac {}{}
-(I_{\alp}(x) -\alp ) (I_{\alp}(y) -\alp )\ha + (I_{\be}(x) -\be
)(I_{\be}(y) -\be ) \hb \rb,
$$
where
$$
E(g_N(X_1))=0, \qquad E(\psi_N(X_1,X_2))=0,
$$
$$
E(\psi_N(X_1,X_2)|X_2)=0   \quad a.s.
$$
The local smoothness assumption of our theorem directly yields that the
distribution of r.v. $g_N(X_1)=
\frac 1{\sqrt{N}}(W_1-\mu_W)$ has a nontrivial absolutely continuous
component and Cram\'er's condition
$$
(C) \qquad \limsup_{|t|\to \infty} |E\exp \{ it \sqrt{N}
g_N(X_1)\}|< 1
$$
is satisfied. Since the functions $\sqrt{N}g_N(x)$ and $N^{3/2}\psi_N(x,y)$
are both bounded, we trivially have that
$$
\be_4 =E\lr \sqrt{N} g_N(X_1)\rr^4 < \infty ,
$$
$$
\qquad \gamma_3 = E \lv N^{3/2}\psi_N (X_1,X_2) \rv^3 < \infty .
$$
Therefore, we can apply Thm.1.2 of Bentkus, G\" otze and van Zwet
\cite{begz} (note that the quantity $\Delta^2_3$ appearing in
Thm.1.2 of Bentkus et. all \cite{begz} is zero in our case).
Define $F_N(x)= \Phi (x) -\phi (x) \frac {\lambda_1
+3\lambda_2}{6\sqrt{N}}(x^2-1)$, where $\lambda_1$ and $\lambda_2$
as in (2.5) (cf.also (4.10)). Then by Thm.1.2 (Bentkus et.all
\cite{begz})
$$
\sup_{x\in R} \lv P\lf \frac {L_N+U_N}{\si_W} \le x\rf - F_N(x) \rv = O(N^{-1}).
$$
For $R_N$ we have the bound (4.12), that is $|R_N|=O((\log
N)^{5/4} N^{-3/4})$ with probability $1-o(N^{-c})$ for every
$c>0$. Therefore, as $F'_N(x)$ and $xF'_N(x)$ are bounded
functions, we obtain on the r.h.s. of (4.19)
$$
F_N(x) - \frac {\sqrt {N}\be_N}{\si_W} \phi (x) + O((\log N)^{5/4}
N^{-3/4})
$$
$$
= G_N (x) +O((\log N)^{5/4} N^{-3/4}).
$$
This proves (2.12) and Theorem 2.1.  $\Box$
\subsection*{5. Proof of Theorem 2.2}
%{\bf 4. Proof of Theorem 1.2.}
\ \ \ \ \ \ Let $S_N^2$ be (cf.(2.14)) the plug in  estimator for
$\si_W^2$ (cf.(2.3)). The following lemma is a modification of
Lemma 4.3 of Putter and van Zwet \cite{putz}, appropriate for our
purposes.

\bigskip
{\sc Lemma} 5.1.\ \ {\it Suppose that the assumptions of Theorem
2.1 are satisfied. Then
$$
P\lr |S_N^2-\si_W^2-V_N|> A(\lon )^{3/4}\rr =O(N^{-c})\leqno (5.1)
$$
as $\nfty$ for every $c>0$ and some $A>0$, not depending on $N$, where
$$
V_N=V_{N,1}+V_{N,2}\ , \leqno (5.2)
$$
$$
V_{N,1}= 2\alp \ha \frac {\na-\alp N}N [\mu_W-\Ha]+
2(1-\be) \hb \frac {\nb-\be N}N [\mu_W-\Hb],
$$
$$
V_{N,2}=\frac 1N \sum_{i=1}^N\lb (\wi -\mu_W)^2-\si_W^2\rb\ .
$$
Moreover,
$$
E(V_N)=0\ ; \ \ \ E(V_N^2)=O(N^{-1})  \leqno (5.3)
$$
as $\nfty$.}

This lemma essentially asserts that the difference between
$\si_W^2$ and its estimator $S_N^2$ can be
expressed as a sum of i.i.d. r.v.'s plus a remainder term which is
of negligible order for our purposes.

\bigskip
{\sc Proof}. Define the auxiliary quantity
$$
S_W^2 = \frac 1N \sum_{i=1}^N W_i^2 - \lr \frac 1N \sum_{i=1}^N
W_i \rr^2
$$
$$
=\frac {\na}N \Had +\frac 1N \sum_{i=\na +1}^{\nb}\xind + \frac
{N-\nb}N \Hbd\ - \lr \frac {\na}N \Ha + \frac 1N \sum_{i=\na
+1}^{\nb} \xin + \frac {N-\nb}N \Hb \rr^2 .
$$
First we prove that
$$
S_N^2=S_W^2+V_{N,1}+R_{N,1}, \leqno (5.4)
$$
Here and elsewhere  $R_{N,1}$, $R_{N,1}^{(r)}$, $r=1,2,\dots$
denote the remainder terms of Bahadur's order, satisfying (3.2).
We have
$$
S_N^2-S_W^2= \lb \frac kN \xknd +\frac 1N \sum_{i=k+1}^{m-1} \xind
+ \frac {N-m+1}N \xmnd  \rp \qquad \qquad \qquad \qquad \leqno
(5.5)
$$
$$
\qquad \qquad \qquad \lp \qquad - \frac {\na}N \Had - \frac 1N
\sum_{i=\na+1}^ {\nb} \xind -\frac {N-\nb}N \Hbd \rb + \lb \lr
\frac 1N \sum_{i=1}^N W_i \rr^2- \hat \mu_W^2 \rb .
$$
Rewrite the term within the first square brackets on the r.h.s. of
(5.5) as
$$
\frac kN (\xknd -\Had) +sign(\na-k)\frac 1N \sum_{i=(k\wedge
\na)+1}^{k\vee \na} (\xind -\Had)
$$
$$
+ \frac {N-m+1}N (\xmnd -\Hbd) - sign(\nb-m+1) \frac 1N
\sum_{i=m\wedge (\nb +1)}^{(m-1)\vee \nb} (\xind -\Hbd)
$$
with $sign(0)=0$ (cf. the proof of Lemma 4.1). By Lemmas 3.1 and
3.2 this expression is equal to
$$
-2\alp \Ha \ha \frac {\na-\alp N}N- \frac {(\na -\alp
N)^2}{N^2}\Ha \ha +R_{N,1}^{(1)} - \leqno (5.6)
$$
$$
- 2(1-\be)\Hb \hb \frac {\nb-\be N}N + \frac {(\nb -\be
N)^2}{N^2}\Hb \hb +R_{N,1}^{(2)},
$$
and by Bernstein's inequality for the binomial r.v.'s $\na$ and
$\nb$ the latter formula reduces to
$$
-2\alp\Ha\ha\frac{\na-\alp N}N - 2(1-\be)\Hb\hb\frac{\nb-\be N}N +
R_{N,1}^{(3)}. \leqno (5.7)
$$
Now we consider the term within the second square brackets on the
r.h.s. of (5.5). Arguing as before, we can rewrite this expression
as
$$
\ \ \lr \frac 2N\sum_{i=1}^N \wi -\alp\ha\frac {\na-\alp N}N-
(1-\be)\hb\frac {\nb-\be N}N +  R_{N,1}^{(4)}\rr \leqno (5.8)
$$
$$
\cdot \lr \alp\ha\frac {\na-\alp N}N+ (1-\be)\hb\frac {\nb-\be N}N
+ R_{N,1}^{(5)}\rr  \qquad \qquad
$$
$$
= \frac 2N \lr\sum_{i=1}^N \wi\rr \lr \alp\ha\frac {\na-\alp N}N+
(1-\be)\hb\frac {\nb-\be N}N \rr  +  R_{N,1}^{(6)}\ .
$$
The relations (5.6)--(5.8) together imply that
$$
S_N^2-S_W^2= V_{N,1}+R_N +R_{N,1}^{(7)}, \leqno (5.9)
$$
where
$$
R_N= 2\lb \alp \ha \frac {\na-\alp N}N +
(1-\be) \hb \frac {\nb-\be N}N \rb \frac 1N \sum_{i=1}^N
(\wi -\mu_W)\ .
$$
Note that the $\wi$, $i=1,\dots ,N$, are bounded i.i.d. $r.v.$'s.
Therefore by Hoeffding's inequality $\frac 1N \lv \sum_{i=1}^N
(\wi -\mu_W)\rv =O\lr (\lon)^{1/2}\rr$ as $\nfty$ with probability
$1-o(N^{-c})$ for every $c>0$. Combining the latter bound with
Bernstein's inequality for the binomial r.v.'s $\na$ and $\nb$, we
obtain that $|R_N|=O(\lon)$ with probability $1-o(N^{-c})$ for
every $c>0$. Therefore (5.9) implies (5.4).

Next we prove that
$$
S_W^2=\si_W^2+V_{N,2}+R_{N,2}, \leqno (5.10)
$$
where $|R_{N,2}|=O(\lon)$ with probability
$1-o(N^{-c})$ for every $c>0$. We have
$$
S_W^2-\si_W^2-V_{N,2}=S_W^2- \frac 1N \sum_{i=1}^N (\wi -\mu_W)^2=
-(\overline W-\mu_W)^2 =R_{N,2}.
$$
An application of Hoeffding's inequality to the bounded i.i.d.
$r.v.$'s  $\wi$ (cf. \cite{hoef})   proves (5.10). Relations (5.4)
and (5.10) together imply (5.1). The lemma is proved.  $\Box$

\bigskip
Now we turn to the proof of our result concerning the Studentized
version of trimmed mean.

\bigskip
{\sc Proof of Theorem} 2.2. \ \ Our proof of this theorem closely
resembles the proof of Theorem 1.2 of Putter and van Zwet
\cite{putz}. For the $df$ $F_{N,S}(x)$ (cf.(2.15)) of a
Studentized trimmed mean we have
$$
F_{N,S}(x) = P\lf \frac {L_N+U_N}{S_N} \le  \lr 1+O(N^{-1})\rr \lb
x -\frac {N^{1/2}\be_N + O(N^{-1})}{S_N} \rb +\frac {R_{N,1}}{S_N}
\rf  \leqno (5.11)
$$
(cf.(4.19)). Here and elsewhere $R_{N,1}$ denotes a remainder,
which satisfies (4.12) and which can be different from line to
line. Lemma 5.1 and Hoeffding's inequality for r.v. $V_N$ together
imply that $\lv \frac 1{S_N}-\frac 1{\si_W}\rv =O((\lon)^{1/2})$
with probability $1-O(N^{-c})$ as $\nfty$ for every $c>0$ (cf.also
Lemma 6.2, below). Therefore, the r.h.s. of (5.11) equals to
$$
P\lf \frac {L_N+U_N}{S_N} \le  \lr 1+O(N^{-1})\rr \lb x -\frac
{N^{1/2} \be_N }{\si_W} \rb + R_{N,1}  \rf .  \leqno (5.12)
$$
Our aim now is to prove that
$$
\sup_{x \in R}|F_{N,S}(x) - H_N(x)| = O\lr(\log
N)^{5/4}/N^{3/4}\rr \leqno (5.13)
$$
as $\nfty$ (this implies (2.17)). Define $\tilde H_N(x) = H_N(x) +
\si^{-1}_W \sqrt {N} \be_N \phi (x)$ (i.e. $\tilde H_N(x)$ is
$H_N(x)$ without bias term). Since $H'_N(x)$ and $xH'_N(x)$ are
bounded, relations (5.11) and (5.12) imply that it is sufficient
to show that
$$
\sup_{x \in R}|F_{(L_N+U_N)/S_N}(x) - \tilde H_N(x)| = O\lr(\log
N)^{5/4}/N^{3/4}\rr, \leqno (5.14)
$$
where $F_{(L_N+U_N)/S_N}(x)=P\lr (L_N+U_N)/S_N \le  x \rr$. An
application of the Lemma 5.1 yields that
$$
F_{(L_N+U_N)/S_N}(x) = P\lr \frac {L_N+U_N}{\si_W} \le  x
\frac {(\si^2_W+V_N+R_N)^{1/2}}{\si_W}\rr,
$$
where $R_N$ is a remainder of Bahadur's order (i.e. satisfying
(3.2)). Since $x\tilde H'_N(x)$ is bounded, it is sufficient to
prove (5.14) with $F_{(L_N+U_N)/S_N}(x)$ replaced by
$$
P\lr \frac {L_N+U_N}{\si_W} \le x \frac {(\si^2_W+V_N )^{1/2}}
{\si_W} \rr =  P\lr \frac {L_N+U_N}{\si_W} - x \lf \lr 1+
\frac {V_N}{\si_W^2}\rr^{1/2} - 1 \rf \le x \rr.
$$
Following Putter and van Zwet \cite{putz}, we also use the
inequality $1+\frac z2 -\frac{z^2}4 \le (1+z)^{1/2} \le 1 +\frac
z2$ ($|z|\le\frac 45$) to find that $\frac {V_N}{2\si_W^2} - \frac
{V_N^2}{4\si_W^4} \le \lr 1+ \frac {V_N}{\si_W^2}\rr^{1/2} - 1 \le
\frac {V_N}{2\si_W^2}$ (with probability $1-O(N^{-c})$, $c>0$).
Since by Hoeffding's inequality $V_N^2=O(\lon)$ with probability
$1-O(N^{-c})$ for every $c>0$, we can replace
$F_{(L_N+U_N)/S_N}(x)$ in (5.14) by $P\lr \frac {L_N+U_N}{\si_W} -
x \frac {V_N}{2\si_W^2} \le x \rr$. Now it remains to show that
$$
\sup_{x\in R}\lv P\lr \frac {L_N+U_N}{\si_W} - x \frac
{V_N}{2\si_W^2} \le x \rr - \tilde H_N(x) \rv = O\lr (\log
N)^{5/4}/N^{3/4} \rr , \leqno (5.15)
$$ as $\nfty$. First we prove
(5.15), taking supremum for $x : |x|<\log N$ (cf. \cite{putz}).
Note that $U_{Nx}= \frac {L_N+U_N}{\si_W} - x \frac
{V_N}{2\si_W^2}$ is a centered $U$-statistic of degree two with
bounded (uniformly for all $x$: $|x|<\log N$) kernel. Moreover,
$U_{Nx}$ has a nontrivial absolutely continuous component and
Cram\'er's condition is satisfied. Theorem 1.1 of Bentkus, G\"otze
and van Zwet \cite{begz} now yields that
$$
\sup_{|x| < \log N}\lv P\lr \frac {L_N+U_N}{\si_W} - \frac
{xV_N}{2\si_W^2} \le x \rr - \tilde G_N(x) \rv = O\lr N^{-1} \rr ,
\leqno (5.16)
$$
where $ \tilde G_N(x) = \Phi \lr\frac x{\si_x}\rr - \frac
{k_{3x}}{6 \si_x^3 } \lb \lr \frac x{\si_x}\rr^2 - 1 \rb \phi \lr
\frac x{\si_x}\rr$ with $\si^2_x = Var(U_{Nx}) = E \lr \frac
{L_N+U_N}{\si_W} - \frac {xV_N}{2\si_W^2} \rr^2$ and $k_{3x} = E
\lr \frac {L_N+U_N}{\si_W} - \frac {xV_N}{2\si_W^2}\rr^3$. Using
the formulas (4.3)--(4.5) and the relations (5.2)--(5.3), we find
that $\si_x^2 = 1+O\lr \frac {\log N}{\sqrt{N}} \rr$ and $k_{3x} =
\frac {\lambda_1 +3\lambda_2}{\sqrt{N}} + O\lr \frac {\log N}{N}
\rr$. Therefore
$$
\tilde G_N(x) = \Phi \lr\frac x{\si_x}\rr - \frac {\lambda_1
+3\lambda_2} {6\sqrt{N}}(x^2-1)\phi (x) +O\lr \frac {\log N}{N}
\rr \leqno (5.17)
$$
(for $|x|<\log N$), that is $\si_x$ influences the form of EE only
through the term $\Phi\lr\frac {x}{\si_x}\rr$ (cf. \cite{putz}).
For $\si^2_x$ we can write $\si_x^2 =E \lr \frac {L_N+U_N}{\si_W}
- \frac {xV_N}{2\si_W^2}\rr^2 = 1- x\si_W^{-3} E[(L_N +U_N)V_N] +
O\lr \frac {\log^2 N}{N} \rr$. As $U_N$ and $V_N$ are
uncorrelated, using formulas (4.3)--(4.4) and (5.2), we can write
$E[(L_N +U_N)V_N] = E(L_N V_N) = \frac 1{\sqrt{N}} (\gamma_{3,W} +
2\delta_{2,W})$. Thus, we obtain that $\si_x^2 = 1- \frac
{x(\lambda_1 + 2\lambda_2)}{\sqrt{N}} + O\lr \frac {\log^2 N}{N}
\rr$ (cf. notations (2.3)--(2.5)). This implies that
$$
\Phi\lr \frac {x}{\si_x}\rr = \Phi (x)  + \phi (x) \frac 12 \frac
{x^2(\lambda_1+2\lambda_2)}{\sqrt{N}}+ O\lr \frac {\log^2 N}N\rr.
\leqno (5.18)
$$
Relations (5.17) and (5.18) together yield that $\tilde G_N(x) =
\tilde H_N(x) + O \lr \frac {\log^2 N}{N} \rr$ for $|x|<\log N$.
To treat the case $|x|\ge \log N$, we use the same arguments as in
 [16, p. 1561] to find that $\sup_{x\in R}\lv P\lr \frac
{L_N+U_N}{\si_W} - \frac {xV_N}{2\si_W^2} \le x \rr - \tilde
H_N(x) \rv = O\lr \frac {\log^2 N}N \rr$. This proves (5.15) and
the theorem. $\Box$
\subsection*{6. Proof of Theorem 2.3}
%{\bf 5. Proof of Theorem 1.3.}
\ \ \ \  In this section we state and prove two lemmas on the
consistency of the estimators for $\lambda_1$, $\lambda_2$ and
$\be_N$. The validity of Theorem 2.3 follows directly from
Theorems 2.1, 2.2 and these lemmas. In the first lemma we obtain
the rate of convergence for our kernel estimates of the density
evaluated at given quantiles, defined by (2.18).

\bigskip
{\sc Lemma} 6.1. {\it Suppose that $f=F'$ exists in a neighborhood
of $\xia $ and satisfies a Lipschitz condition. In addition we
assume that $f(\xia)>0$. Then
$$
P\lr |\hat f(\xia) - f(\xia)|>A(\log N)^{1/2}/N^{1/4}\rr =
O(N^{-c}) \leqno (6.1)
$$
as $\nfty$, for every $c>0$ and some $A>0$, not depending on $N$.
}
\bigskip

{\sc Proof}. Define random quantities
$$
\nu_{k,N}= \sharp \lf X_i : |X_i-\xkn |\le N^{-1/4}/2 \rf,\ \
\nu_{\alp,N}= \sharp \lf X_i : |X_i-\xia |\le N^{-1/4}/2 \rf.
\leqno (6.2)
$$
Note that $E\nu_{\alp,N}= N\int_{\xia - N^{-1/4}/2}^{\xia +
N^{-1/4}/2} f(x) \, dx $, and one can write
$$
\hat f(\xia) - f(\xia) = N^{-3/4}\nu_{k,N}- f(\xia) \qquad \qquad
\qquad \qquad \qquad \quad \ \leqno (6.3)
$$
$$
= N^{-3/4}\nu_{\alp,N}+ N^{-3/4}(\nu_{k,N}-\nu_{\alp,N}) - f(\xia)
= Q_{1,N} +Q_{2,N}+Q_{3,N},
$$
where
$$
Q_{1,N}= N^{-3/4}(\nu_{\alp,N}- E\nu_{\alp,N}), \qquad
Q_{2,N}= N^{-3/4}(\nu_{k,N}-\nu_{\alp,N}),
$$
$$
Q_{3,N}=N^{1/4}\int_{\xia - N^{-1/4}/2}^{\xia + N^{-1/4}/2}
(f(x) - f(\xia) \, dx .
$$
For $Q_{1,N}$ we can write $Q_{1,N}= N^{1/4}\lr \bar \nu_{\alp,N}-
E\bar \nu_{\alp,N}\rr$, where $\bar \nu_{\alp,N} =\frac 1N
\sum_{i=1}^N I_{\{2N^{1/4}|X_i-\xia|\le 1\} }$ is a mean of i.i.d.
bounded r.v.'s. Therefore, by Hoeffding's inequality
$$
P\lr |Q_{1,N}|>A_1 (\log N)^{1/2}/N^{1/4} \rr = O(N^{-c}) \leqno
(6.4)
$$
for every $c>0$, as $\nfty$. Here and elsewhere $A_i$, $i=1,2,\dots$ denote
positive constants, not depending on $N$. Since $P(|\xkn -\xia|>A_2
(\lon)^{1/2} )=O(N^{-c})$, for $Q_{2,N}$ we have with probability
$1-O(N^{-c})$
$$
|Q_{2,N}|\le N^{-3/4} (\nu_{l,N}+\nu_{r,N}), \leqno (6.5)
$$
where $\nu_{l,N}= \sharp \lf X_i : |X_i-\xia +N^{-1/4}/2| \rp$ $\lp
\le A_2 (\lon )^{1/2} \rf$, $\nu_{r,N}= \sharp \lf X_i : |X_i-\xia -
N^{-1/4}/2|\rp $ $\lp \le A_2 (\lon )^{1/2} \rf$. Since $(\nu_{l,N}+\nu_{r,N})$ is
a Binomial r.v. with parameter $p_N =O\lr (\lon )^{1/2}\rr$ and
$E(\nu_{l,N}+\nu_{r,N})= O\lr N^{1/2}(\log N)^{1/2}\rr$,
$\si_{\nu_{l,N}+\nu_{r,N}}= O\lr N^{1/4}(\log N)^{1/4}\rr$, by Bernstein
inequality, with probability $1-O(N^{-c})$, we have the following bound
$$
|Q_{2,N}|\le A_3 N^{-1/4} (\log N)^{1/2}. \leqno (6.6)
$$
Finally for $Q_{3,N}$ the Lipschitz condition directly yields that
$$
|Q_{3,N}| \le C N^{1/4} \int_{\xia - N^{-1/4}/2}^{\xia +
N^{-1/4}/2} |x - \xia|\, dx = \frac 14 C N^{-1/4}, \leqno (6.7)
$$
where $C$ is the Lipschitz constant. Relations (6.3)--(6.7) imply
(6.1). The lemma is proved. $\Box$
\bigskip

Let $\mu_{r,W}=EW_i^r=\int_0^1 Q^r(u)\, du$ denotes the $r$-th
moment of $W_i$ for any positive integer $r$ and let $\hat
\mu_{r,W}=\frac kN \xkn^r+ \frac 1N \sum_{i=k+1}^{m-1} \xin^r
+\frac {N-m+1}N \xmn^r$ be the plug in estimator for $\mu_{r,W}$.

\bigskip
{\sc Lemma} 6.2. {\it Suppose that $f=F'$ exists in neighborhoods
of $\xia$ and $\xib$ and satisfies a Lipschitz condition. In
addition we assume that $f(\xi_{\nu})>0$, $\nu=\alp ,\be $. Then
$$
P\lr |\hat \mu_{r,W} - \mu_{r,W}|>A(\lon )^{1/2}\rr = O(N^{-c})
\leqno (6.8)
$$
as $\nfty$ for every $c>0$ with some $A>0$, not depending on $N$.
}

\bigskip

{\sc Proof}. Put $\bar W_r=\frac 1N \sum_{i=1}^N W_i^r$, where
$W_i$ is defined by (2.2), and note that similarly when proving of
lemma 5.1 we can write
$$
\bar W_r=\frac {\na}N \Har + \frac 1N \sum_{i=\na +1}^{\nb} \xin^r
+\frac {N-\nb}N \Hbr.
$$
We have
$$
\hat \mu_{r,W} - \mu_{r,W} = (\hat \mu_{r,W} -\bar W_r) + (\bar
W_r - \mu_{r,W}). \leqno (6.9)
$$
Note that $E\bar W_r = \mu_{r,W}$, therefore by Hoeffding
inequality for the average of i.i.d. bounded r.v.'s we have $
|\bar W_r - \mu_{r,W}| = O\lr (\lon)^{1/2}\rr $ with probability
$1-O(N^{-c})$ for every $c>0$. For $(\hat \mu_{r,W} -\bar W_r)$ on
the r.h.s. of (6.9) we have
$$
\hat \mu_{r,W} -\bar W_r = \frac kN (\xkn^r - \Har) +
sign(\na-k)\frac 1N \sum_{i=(k\wedge \na)+1}^{k\vee
\na}(\xin^r-\Har) \qquad \qquad \qquad
$$
$$
\qquad \qquad \qquad + \frac {N-m+1}N (\xmn^r - \Hbr) - sign(\nb -
m+1)\frac 1N \sum_{i=m\wedge (\nb+1)}^{(m-1)\vee
\nb}(\xin^r-\Hbr).
$$
(cf.(5.6)). By Lemmas 3.1 and 3.2 the last expression equals to
$$
-\alp r \Harl \frac {\na - \alp N}N \ha - \frac {(\na - \alp N)^2}
{2N^2} r\Harl\ha \leqno (6.10)
$$
$$
- (1-\be)r \Hbrl \frac {\nb - \be N}N \hb + \frac {(\nb - \be
N)^2} {2N^2} r\Hbrl\hb + R_N,
$$
where $R_N$ is a remainder term of the Bahadur's order (cf.
(3.2)). Thus, by Bernstein inequality we find that
$$
|\hat \mu_{r,W} -\bar W_r|=O\lr (\lon)^{1/2}\rr \leqno (6.11)
$$
with probability $1-O(N^{-c})$ for every $c>0$. Relations
(6.9)--(6.11) together imply (6.8). The lemma is proved. $\Box$

\subsection*{Appendix}
%\begin{center}
%APPENDIX
%\end{center}
\ \ \ \  In this appendix we first establish an asymptotic
approximation for the bias of $T'_N$ (cf. (2.9)) in estimating of
$\mu (\alp ,\be )$. Secondly we prove that our Theorem 2.2 can not
be inferred from Theorem 1.2 of Putter and van Zwet \cite{putz}
for Studentized symmetric statistics..

\bigskip
{\sc Lemma} A.1. \ \ {\it Suppose the conditions of Theorem 2.1
are satisfied. Then
$$
b_N = \be_N +O(N^{-3/2}), \leqno (A.1)
$$
with $b_N$ and $\be_N$ as in (2.8) and (2.9). }

\bigskip
{\sc Proof}. To begin with we note that $b_N$ (cf. (2.9) can be
written as $B_1 +B_2$ where $B_1 =(\be -\alp) ET'_N - E \lr \frac
1N \sum_{i=[\alp N]+1}^{[\be N]} \xin' \rr$ and $B_2=E \lr \frac
1N \sum_{i=[\alp N]+1}^{[\be N]} \xin' \rr - (\be -\alp)\mu (\alp
, \be)$. First we consider $B_2$. By a simple conditioning
argument we have that $B_2$ equals (with $k=[\alp N]+1$, $m=[\be
N]$)
$$
\frac 1{N} E\lr F^{-1}(\ukn)+F^{-1}(\umn) +(m-k-1) \frac
{\int_{\ukn}^{\umn} F^{-1}(u)\, du}{\umn-\ukn}\rr - (\be -\alp)\mu
(\alp , \be). \leqno (A.2)
$$
Define
$$
I(v_1,v_2)= \frac {\int_{v_1}^{v_2}F^{-1}(u)\, du}{v_1-v_2},
\qquad I(\alp,\be)= \mu (\alp , \be).
$$
The first and second partial derivatives are given by
$$
\lp \frac {\dd I}{\dd v_1} \rv_{(\alp,\be)} =
\frac {-\Ha +\mu (\alp,\be)}{\be-\alp} , \qquad
\lp \frac {\dd I}{\dd v_2} \rv_{(\alp,\be)} =
\frac {\Hb -\mu (\alp,\be)}{\be-\alp}\ ,
$$
$$
\lp \frac {\dd^2 I}{\dd v_1^2} \rv_{(\alp,\be)} = -\frac {2}
{\be-\alp}\lb \frac 1{2f(\xia)}-\frac {\mu (\alp , \be) - \Ha}
{\be-\alp} \rb , \qquad \qquad \quad
$$
$$
\lp \frac {\dd^2 I}{\dd v_2^2} \rv_{(\alp,\be)} = \frac {2}
{\be-\alp}\lb \frac 1{2f(\xib)}+\frac {\mu (\alp , \be) - \Hb}
{\be-\alp} \rb , \qquad \qquad \qquad
$$
$$
\lp \frac {\dd^2 I}{\dd v_1\dd v_2} \rv_{(\alp,\be)} = \frac {\Ha +\Hb -
2\mu (\alp , \be)}{(\be-\alp)^2}. \qquad \qquad \qquad \qquad \qquad
\qquad \quad
$$
A Taylor expansion argument now yields that (A.2) reduces to
$$
\frac 1N E \lr F^{-1}(\ukn)+ \frac {}{}F^{-1}(\umn) \rr +\frac
{m-k-1}{N} \lf \frac {}{} \mu (\alp , \be ) \rp \qquad \qquad
\quad
$$
$$
+ \frac {-\Ha +\mu (\alp,\be)}{\be-\alp} \lr \frac k{N+1} -\alp
\rr + \frac {\Hb -\mu (\alp,\be)}{\be-\alp} \lr \frac m{N+1}-\be
\rr \qquad \qquad \quad
$$
$$
- \frac {1} {\be-\alp}\lb \frac 1{2f(\xia)}-\frac {\mu (\alp ,
\be) - \Ha} {\be-\alp} \rb   \frac {\frac k {N+1}(1-\frac k
{N+1})}{N+2}   \qquad \qquad \qquad \qquad \qquad
$$
$$
+ \frac {1} {\be-\alp}\lb \frac 1{2f(\xib)}+\frac {\mu (\alp ,
\be) - \Hb} {\be-\alp} \rb  \frac {\frac m {N+1}(1-\frac m
{N+1})}{N+2}  \qquad \qquad \qquad \qquad \qquad
$$
$$
+ \lp \lb \frac {\Ha +\Hb - 2\mu (\alp , \be)}{(\be-\alp)^2}\rb
\frac {\frac k{N+1}(1-\frac m{N+1})}{N+2} + O(N^{-3/2}) \rf
-(\be-\alp) \mu (\alp,\be),
$$
which easily leads to
$$
\ \frac 1N \lf \frac {}{}\Ha (\alp N -[\alp N]) - \Hb (\be N -[\be
N])  \rp  \qquad \qquad  \leqno (A.3)
$$
$$
- \lp\frac 1{2f(\xia)} \alp(1-\alp) + \frac 1{2f(\xib)} \be
(1-\be) \rf +O(N^{-3/2}).
$$
For $B_1$ we have
$$
B_1 = \frac {(\be N - [\be N]) - (\alp N - [\alp N])}{[\be N] -
[\alp N]} E \lr \frac 1N \sum_{i=k}^m \xin' \rr \qquad \qquad
\qquad \qquad \qquad
$$
$$
= \frac {(\be N - [\be N]) - (\alp N - [\alp N])}{[\be N] - [\alp
N]} \lr \frac {}{}(\be-\alp) \mu (\alp,\be) +\be_N +O(N^{-3/2})
\rr
$$
$$
= \frac {1}{N} \lr \frac {}{} (\be N - [\be N]) - (\alp N - [\alp
N]) \rr \mu (\alp,\be) + O(N^{-2}). \quad \qquad \qquad \ \
$$
This together with (A.2)--(A.3) implies (A.1). The lemma is proved. $\Box$

\bigskip
Consider a trimmed mean $T_N$ as in (4.1). Let $T_{N\Omega_k}$ is
defined as in (1.8) of Putter and van Zwet \cite{putz}. We prove
the following assertion.

\bigskip
{\sc Lemma A.2.} {\it Suppose that the conditions of Theorem 2.1
hold. Then
$$
\sum_{k=3}^N {N-2 \choose k-2}ET^2_{N\Omega_k} = N^{-3} \lr \frac
{\alp^2(1-\alp)^2}{f^2(\xia)}+ \frac
{\be^2(1-\be)^2}{f^2(\xib)}\rr + o(N^{-3}) \leqno (A.4)
$$
as $\nfty$.
}

\bigskip
Relation (A.4) directly yields that in the second condition of
(1.18) in Theorem 1.2 of Putter and van Zwet \cite{putz} is not
satisfied for a Studentized trimmed mean, as Putter and van Zwet
\cite{putz} require that the l.h.s. of (A.4) is of order
$N^{-7/2}$, instead of $N^{-3}$ as in our relation (A.4).

\bigskip
{\sc Proof}. In Putter's Ph.D thesis \cite{putt} it was proved
that if $T_N$ is a linear combination of order statistics, then
$$
\sum_{k=3}^N {N-2 \choose k-2}ET^2_{N\Omega_k} = E(Z_N -
E(Z_N|U_{N-1},U_N))^2 \leqno (A.5)
$$
$$
\qquad \qquad \qquad = EZ^2_N - E(T_{N,(1,2)})^2, \
$$
(cf. (3.5.17), Putter \cite{putt}), where $T_{N\Omega_k}$,
$T_{N,(1,2)}$ are defined as in (1.8) of Putter and van Zwet
\cite{putz}, $Z_N$ is a r.v. defined as in (4.21) of van Zwet
\cite{zwet}, $U_1,\dots ,U_N$ are uniformly on (0,1) distributed
r.v.'s. Let $R_j$ denotes the rank of $U_j$ among $U_1,\dots
,U_N$, \ $K_1 = R_{N-1}\wedge R_N$, \ $K_2 = R_{N-1}\vee R_N$.
Take $X_{0:N}=-\infty$, $X_{N+1:N}=+\infty$ (cf. van Zwet
\cite{zwet}). Let the functions $G$, $H$, $M$ are defined as in
(4.17) of van Zwet \cite{zwet}, and define in addition the
functions $G_1$ and $H_1$ by $G_1(x)=\int_{-\infty}^x F^2(y) \,
dy$, $H_1(x)=\int_x^{\infty} (1- F(y))^2 \, dy$. Then formula
(4.21) of van Zwet \cite{zwet} reduces to
$$
N^{1/2}Z_N=-\sum_{j=1}^{K_1} (c_{j+1}-c_j)(G_1(X_{j:N})-
G_1(X_{j-1:N}))
$$
$$
+ \sum_{j=K_1}^{K_2-1} (c_{j+1}-c_j)(M(X_{j+1:N})- M(X_{j:N})) -
\sum_{j=K_2}^N (c_j-c_{j-1})(H_1(X_{j:N})- H_1(X_{j+1:N}))
$$
(cf. Gribkova \cite{gri}), where in the trimmed mean case ($c_j=1$
for $k\le j \le m$ and $c_j=0$ for $j<k$, $j>m$) there are only
two nonzero summands, which depend on $K_1$ and $K_2$ (the event
\{$K_1=k-1$ or $K_2=m+1$\} is negligible for our aims because its
probability is $O(N^{-1})$, cf. below). For instance, when $K_2<k$
(which happens with probability $P(K_2<k)=\alpha^2+O(N^{-1}))$,
the value of $N^{1/2}Z_N$ equals
$$ -[H_1(X_{k:N})-H_1(X_{k+1:N})]+[H_1(X_{m+1:N})-H_1(X_{m+2:N})]
\stackrel{d}{=}
$$
$$-[H_1\circ F^{-1}(U_{k:N})-H_1\circ F^{-1}(U_{k+1:N})]+
[H_1\circ F^{-1}(U_{m+1:N}) -H_1\circ F^{-1}(U_{m+2:N})],$$ where
$U_{i:N}$ are order statistics of r.v.'s $U_i,~i=1,\ldots,N$.
Application of a two term Taylor expansion of the function
$H_1\circ F^{-1}$ in neighborhoods of $\alpha$ and $\beta$
respectively, together with the well-known facts that
$E(s_i^2)=\frac{2}{(N+2)(N+1)}$,
$E(s_is_j)=\frac{1}{(N+2)(N+1)}(i\neq j)$, where
$s_i=U_{i:N}-U_{i-1:N},~i=1,\ldots,N+1$, $U_{0:N}=0$,
$U_{N+1:N}=1$, yields that $E(Z_N^2 | K_2<k)= \frac{2}{N^3}\left(
\frac{(1-\alpha)^4}{f^2(\alpha)}+ \frac{(1-\beta)^4}{f^2(\beta)}-
\frac{(1-\alpha)^2(1-\beta)^2}{f(\alpha)f(\beta)}\right)+o(N^{-3})$,
where $P(K_2<k)=\alpha^2+O(1/N)$. Analyzing in similar fashion the
other possibilities for $K_1$ and $K_2$, we find that
$$
EZ^2_N = \frac 2{N^3} \lr \frac {\alp^2(1-\alp)^2}{f^2(\xia)} -
\frac {\alp^2(1-\be)^2}{f(\xia)f(\xib)} + \frac
{\be^2(1-\be)^2}{f^2(\xib)}\rr + o(N^{-3}), \leqno (A.6)
$$
as $\nfty$. Next we consider $T_{N,(1,2)}$. By formula (2.11) of
Putter and van Zwet \cite{putz} we have
$$
N^{1/2} T_{N,(1,2)} = -\int_0^1 (I_{[U_1,1)}(t)-t)
(I_{[U_2,1)}(t)-t){N-2 \choose k-2}t^{k-2}(1-t)^{N-k} \,
dF^{-1}(t)
$$
$$
\qquad \qquad \qquad \quad \ \  + \int_0^1 (I_{[U_1,1)}(t)-t)
(I_{[U_2,1)}(t)-t){N-2 \choose m-1}t^{m-1}(1-t)^{N-m-1} \,
dF^{-1}(t).
$$
Define $\Delta F_{i,N}(x) = F_{i-1,N}(x)- F_{i,N}(x)$, where
$F_{i,N}(x)= P(\xin \le x)$. The latter relation implies that
$E(T_{N,(1,2)})^2$ equals to
$$
\frac 2N \int_{-\infty}^{\infty} \int_{-\infty}^z\lb
-\int_{-\infty}^{\infty} (I_{[y,+\infty)}(x)-F(x))
(I_{[z,+\infty)}(x)-F(x))\Delta F_{k-1,N-2}(x) \, dx \  \rp \leqno
(A.7)
$$
$$
+ \lp \int_{-\infty}^{\infty} (I_{[y,+\infty)}(x)-F(x))
(I_{[z,+\infty)}(x)-F(x))\Delta F_{m,N-2}(x) \, dx  \rb^2 \, dF(y)
\, dF(z)
$$
$$
\quad \ = \frac 2N \int_{-\infty}^{\infty} \int_{-\infty}^z \lb
I_{k-1}(y,z) \rb^2 \, dF(y)\, dF(z) + \frac 2N
\int_{-\infty}^{\infty} \int_{-\infty}^z \lb I_m(y,z) \rb^2 \,
dF(y)\, dF(z) \
$$
$$
- \frac 4N \int_{-\infty}^{\infty} \int_{-\infty}^z \lb
I_{k-1}(y,z)I_m(y,z) \rb \, dF(y)\, dF(z),
$$
where $I_{r}(y,z)$, $r=k-1, m$, is defined by
$$
\int_{-\infty}^y \Delta F_{r,N-2}(x) \, dG_1(x) -
\int_y^z \Delta F_{r,N-2}(x) \, dM(x) - \int_z^{\infty}
\Delta F_{r,N-2}(x) \, dH_1(x).
$$
Consider the first term at the r.h.s. of (A.7) (the treatment of
the second and third term is similar). Integrating by parts, we
reduce it to
$$
\frac 2N \int_{-\infty}^{\infty} \int_{-\infty}^z\lb (G_1(y)+M(y))
\Delta F_{k-1,N-2}(y) +(H_1(z)-M(z)) \Delta F_{k-1,N-2}(z)\  \frac
{}{} \rp
$$
$$
- \int_{-\infty}^y G_1(x) \, d(\Delta F_{k-1,N-2}(x)) + \int_y^z
M(x)d(\Delta F_{k-1,N-2}(x))
$$
$$
+ \lp \int_z^{\infty} H_1(x) d(\Delta F_{k-1,N-2}(x)) \rb^2 \,
dF(y)\, dF(z).
$$
Note that the `basic' support of the function $\Delta
F_{k-1,N-2}(x) = F_{k-2,N-2}(x) - F_{k-1,N-2}(x)$ is some interval
$I_{\alp}(A) = [\xia - A(\lon)^{1/2}, \xia + A(\lon)^{1/2}]$ in
the sense that for every \hbox{$c>2$} we have the following bound:
$\sup_{y\in R\setminus I_{\alp}(A)} \Delta F_{k-1,N-2}(y) = O\lr
P\lr |U_{k:N}-\alp | > (\lon )^{1/2}\rr\rr= O(N^{-c})$, where
$A>0$ is some constant, depending only on $c$, $\alp$ and
$f(\xia)$. Moreover, smoothness conditions imply that $\sup_{y\in
I_{\alp}(A)} \Delta F_{k-1,N-2}(y)= O(N^{-1})$ as $\nfty$. Thus,
the last expression reduces to
$$
\frac 2N \int_{-\infty}^{\infty} \int_{-\infty}^z\lb
-\int_{-\infty}^y G_1(x) \, d(\Delta F_{k-1,N-2}(x)) + \int_y^z
M(x)d(\Delta F_{k-1,N-2}(x))  \rp \leqno (A.8)
$$
$$
+ \lp \int_z^{\infty} H_1(x) d(\Delta F_{k-1,N-2}(x)) \rb^2 \,
dF(y)\, dF(z) +o(N^{-3}),
$$
as $\nfty$. Consider the integrand in (A.8) and note that if
$I_{\alp}(A)\subset (-\infty ,y)$, then the integrand equals to
$\lb E(G_1(X_{k-2:N-2})- G_1(X_{k-1:N-2}))\rb^2 +o(N^{-2})= \frac
{\alp^4}{N^2}\frac 1{f^2(\xia)}+o(N^{-2})$, and the corresponding
part of the integral in (A.8) (in the domain where $Y=\min
(X_1,X_2) \ge \xia $) equals to $\frac
{(1-\alp)^2\alp^4}{f^2(\xia)}N^{-3}+o(N^{-3})$. Arguing similarly
for the cases $I_{\alp}(A)\subset (y,z)$ and $I_{\alp}(A)\subset
(z,+\infty)$ (the cases $y\in I_{\alp}(A)$ or $z\in I_{\alp}(A)$
are negligible) we obtain that the quantity (A.8), and hence the
first term at the r.h.s. in (A.7), equals to $\lr \frac {\alp^4
(1-\alp)^2}{f^2(\xia)}+ 2\frac {\alp^2 (1-\alp)^2}{f^2(\xia)}+
\frac {(1-\alp)^4 \alp^2}{f^2(\xia)}\rr N^{-3}+o(N^{-3})= $ $\frac
{\alp^2 (1-\alp)^2}{f^2(\xia)}N^{-3} +o(N^{-3})$. Similarly for
the second term at the r.h.s. of (A.7) we get $\frac {\be^2
(1-\be)^2}{f^2(\xib)}N^{-3} +o(N^{-3})$, and for the third one we
obtain $-2 \frac {\alp^2 (1-\be)^2} {f(\xia)f(\xib)}N^{-3}
+o(N^{-3})$. Together these results give us
$$
E(T_{N,(1,2)})^2 = N^{-3}\lr \frac {\alp^2 (1-\alp)^2}{f^2(\xia)}
-2 \frac {\alp^2 (1-\be)^2} {f(\xia)f(\xib)}+ \frac {\be^2
(1-\be)^2}{f^2(\xib)}\rr + o(N^{-3}) \leqno (A.9)
$$
as $\nfty$. The relations (A.5), (A.6) and (A.9) together imply
(A.4) and the lemma is proved. $\Box$
%\subsection*{Acknowledgement}
%This paper was completed while the first author visited CWI in January-February,
%2002; we thank the Netherlands Organization for Scientific Research (NWO)
%for making this visit to CWI possible.

\end{document}